# АНАЛОГ ЗАДАЧИ ТРИКОМИ ДЛЯ УРАВНЕНИЯ СМЕШАННОГО ТИПА С ДРОБНОЙ ПРОИЗВОДНОЙ РИМАНА-ЛИУВИЛЛЯ

*Окбоев Акмалжон Бахромжонович, PhD,*
*akmaljon12012@gmail.com*
*Институт математики имени В.И.Романовского АН РУз,*
*Ташкент, Узбекистан*

***Аннотация.*** *В работе исследована задача Трикоми для уравнения параболо-гиперболического типа в смешанной области. Параболическая часть рассматриваемого уравнения состоит из дробной производной по Риману-Лиувиллю, а гиперболическая часть состоит из вырождающегося гиперболического уравнения второго рода. Решение поставленной задачи в гиперболической подобласти найдено как решение задачи Коши, а в параболической подобласти - как решение первой краевой задачи. Для доказательства существования решения задачи используется теория интегральных уравнений Вольтерра второго рода.*

***Ключевые слова:*** *Уравнение параболо-гиперболического типа, смешанная область, задача Трикоми, задача Коши, первая краевая задача.*

# AN ANALOG OF THE TRICOMI PROBLEM FOR A MIXED TYPE EQUATION WITH RIEMANN-LIOUVILLE FRACTIONAL DERIVATIVE

*Akmaljon Okboev Bakhromjonovich, PhD,*
*akmaljon12012@gmail.com*
*V.I. Romanovskii Institute of Mathematics, Uzbekistan Academy of Sciences,*
*Tashkent, Uzbekistan*

***Abstract.*** *In this article, the Tricomi problem for a parabolic-hyperbolic type equation in a mixed domain is investigated. Riemann-Liouville fractional derivative participates in the parabolic part of the considered equation, and the hyperbolic part consists of a degenerate hyperbolic equation of the second kind. The solution of the problem in the hyperbolic sub-domain is found as a solution to the Cauchy problem, and in a parabolic sub-domain as a solution to the first boundary value problem. For proving the existence of the solution of the problem, the theory of second kind Volterra integral equations is used.*

**Введение.** В этой работе в области $\Omega = \Omega_1 \cup \Omega_2 \cup AB$ для уравнения

$$0 = \begin{cases} L_1(u) \equiv \dfrac{\partial^2 u}{\partial x^2} - D_{0_y}^{\delta} u - \lambda^2 u, & (x,y) \in \Omega_1, \\ L_{\alpha,\lambda}(u) \equiv \dfrac{\partial^2 u}{\partial x^2} + y\dfrac{\partial^2 u}{\partial y^2} + \alpha\dfrac{\partial u}{\partial y} - \lambda^2 u, & (x,y) \in \Omega_2 \end{cases} \quad (1)$$

сформулируем и исследуем задачи Трикоми, где $\Omega_1$ -область, ограниченная при $y > 0$ отрезками $AB, BB_0, A_0B_0, AA_0$ прямых $y = 0, x = 1, y = 1, x = 0$ соответственно, $\Omega_2$ - область, ограниченная при $y < 0$ дугами $AC$, $BC$, $AB$ характеристик $x - 2\sqrt{-y} = 0$, $x + 2\sqrt{-y} = 1$, $y = 0$ уравнения (1), $\delta \in (0,1)$ $\alpha \in (-1/2, 0)$, а $\lambda$ - действительное или чисто мнимое число, $D_{0x}^{\alpha}\varphi(x)$ - интегро-дифференциальный оператор порядка $\alpha$ в смысле Римана-Лиувилля

$$D_{0x}^{\alpha}\varphi(x) = \begin{cases} \dfrac{1}{\Gamma(-\alpha)} \int_0^x \dfrac{\varphi(t) dt}{(x-t)^{\alpha+1}}, & \alpha < 0, \\ \varphi(x), & \alpha = 0, \\ \dfrac{d^n}{dx^n} D_{0x}^{\alpha-n}\varphi(x), & \alpha > 0. \end{cases}$$

**Постановка задачи**

**Определение 1.** *Регулярным в области $\Omega_1$ решением уравнения (1), называется функция $u(x,y)$, удовлетворяющая в области $\Omega_1$ уравнению (1) и следующим условиям $D_{0y}^{\delta-1} u(x,y) \in C(\overline{\Omega_1})$, $u_{xx}(x,y)$, $D_{0y}^{\delta} u(x,y) \in C(\Omega_1)$.*

**Задача $T_0$.** Требуется определить функцию $u(x,y)$, обладающую следующими свойствами: а) $u(x,y)$ является регулярным решением уравнения (1) в области $\Omega_1$ и решением класса $R_{00}^{\lambda}$ в области $\Omega_2$; б) на линии вырождения выполняется условие склеивания

$$\lim_{y \to -0} u(x,y) = \lim_{y \to +0} y^{1-\delta} u(x,y), \quad 0 \leq x \leq 1;$$

$$\lim_{y \to -0} (-y)^{\alpha} \frac{\partial}{\partial y}\left[u(x,y) - A_{\alpha}^{-}(\tau, \lambda)\right] = \lim_{y \to +0} y^{1-\delta} \frac{\partial}{\partial y}\left[y^{1-\delta} u(x,y)\right], \quad 0 < x < 1;$$

в) на границе области $\Omega$ удовлетворяет граничным условиям

$$u|_{AA_0} = \varphi_0(y), \quad u|_{BB_0} = \varphi_1(y), \ 0 \leqslant y \leqslant 1; \tag{2}$$

$$u|_{AC} = \psi(x), \ 0 \leqslant x \leqslant 1/2, \tag{3}$$

где $A_{\alpha}^{-}(\tau, \lambda)$ – определяется формулой

$$A_{\alpha}^{-}(\tau, \lambda) = \gamma_1 \int_0^1 \tau(\zeta) \left[z(1-z)\right]^{\beta} \overline{J}_{\beta}(\sigma) dz + \frac{8\gamma_1 y}{(1+\beta)(1+2\beta)} \times$$

$$\times \int_0^1 \left(\lambda^2 - \frac{d^2}{d\zeta^2}\right) \tau(\zeta) \left[z(1-z)\right]^{1+\beta} \overline{J}_{1+\beta}(\sigma) dz,$$

$\gamma_1 = \Gamma(1+2\alpha)/\Gamma^2(1/2+\alpha)$, $\sigma = 4\lambda\sqrt{-yz(1-z)}$, $\zeta = x - 2\sqrt{-y}(1-2z)$, $\tau(x) = u(x,-0)$, $J_{\gamma}(z)$ - функция Бесселя первого рода, $\overline{J}_{\gamma}(z) = \Gamma(\gamma+1)(z/2)^{-\gamma} J_{\gamma}(z)$, т.е.

$$\bar{J}_\gamma(z) = \Gamma(\gamma+1)\sum_{m=0}^{\infty}\frac{(-1)^m (z/2)^{2m}}{m!\Gamma(m+\gamma+1)}, \ \gamma \neq -1,-2,-3,...., \ \text{а} \ \varphi_0(y), \varphi_1(y) \ \text{и} \ \psi(x) -$$

заданные непрерывные функции.

Отметим, что Н.К.Мамадалиев [1, 2] исследовал различные задачи для уравнения (1) при различных значениях $\alpha$, когда $\delta = 1$, $\lambda = 0$. В работе [3] изучено нелокальная краевая задача для уравнения (1) при $\alpha = \alpha_0 - n, \alpha_0 \in (1/2,1), n = 2,3,...$ и $\delta = 1$. В работе [4] поставленна и исследована задача Трикоми для уравнения (1) при $\alpha = \alpha_0 - n, \alpha_0 \in (1/2,1), n = 2,3,...$ и $\delta = 1$.

**Свойства некоторых операторов с функциями Бесселя в ядрах.**

Рассмотрим следующие интегральные операторы [5]:

$$A_{kx}^{1,\lambda}[g(x)] = g(x) - \int_k^x g(t)\frac{t-k}{x-k}\frac{\partial}{\partial t}J_0\left[\lambda\sqrt{(x-k)(x-t)}\right]dt, \quad (4)$$

$$B_{kx}^{1,\lambda}[g(x)] = g(x) + \int_k^x g(t)\frac{\partial}{\partial x}J_0\left[\lambda\sqrt{(k-t)(x-t)}\right]dt. \quad (5)$$

**Свойство 1.** *Если $g(x) \in C(0,1) \cap L_1[0,1]$, то выражения $A_{kx}^{1,\lambda}[g(x)]$ и $B_{kx}^{1,\lambda}[g(x)]$ будут определены в $(0,1)$ и принадлежат классу $C(0,1)$.*

Справедлива следующая теорема и лемма:

**Теорема 1**[5]**.** *Если $g(x) \in C[0,1]$, то для любых $k \in [0,1]$ и $x \in (0,1)$ справедливы следующие равенства: $A_{kx}^{1,\lambda}B_{kx}^{1,\lambda}[g(x)] = g(x)$, $B_{kx}^{1,\lambda}A_{kx}^{1,\lambda}[g(x)] = g(x)$, т.е. в классе непрерывных на $[0,1]$ функций операторы (4) и (5) являются взаимно обратными.*

**Лемма 1**[5]**.** *При $\beta < 1$ и $x \in [0,1]$ справедливы равенства*

$$\int_0^x (x-t)^{-\beta}\bar{J}_{-\beta}\left[\lambda\sqrt{t(x-t)}\right]g(t)dt = \Gamma(1-\beta)D_{0x}^{\beta-1}\left\{B_{0x}^{1,\lambda i}[g(x)]\right\}, \quad (6)$$

$$\int_x^1 (t-x)^{-\beta}\bar{J}_{-\beta}\left[\lambda\sqrt{(t-x)(1-t)}\right]g(t)dt = \Gamma(1-\beta)D_{x1}^{\beta-1}\left\{B_{x1}^{1,\lambda i}[g(x)]\right\}.$$

**Исследование задачи** $T_0$

Рассмотрим уравнение (1) в области $\Omega_2$, т.е. рассмотрим уравнение $L_{\alpha,\lambda}(u) = 0$. Непрерывное в $\bar{\Omega}_2$ решение видоизмененной задачи Коши для уравнения $L_{\alpha,\lambda}(u) = 0$, с начальными данными

$$u(x,y)\big|_{y=0} = \tau(x),\ 0 \leqslant x \leqslant 1; \qquad (7)$$

$$\lim_{y \to -0} (-y)^\alpha (\partial/\partial y)\left[u - A_\alpha^-(\tau,\lambda)\right] = \nu(x),\ 0 < x < 1, \qquad (8)$$

в характеристических переменных $\xi = x - 2\sqrt{-y}$, $\eta = x + 2\sqrt{-y}$ имеет вид

$$u(x,y) = A_\alpha^-(\tau,\lambda) - 2^{-2+4\beta}\gamma_2 \int_\xi^\eta (\eta-t)^{-\beta}(t-\xi)^{-\beta}\overline{J}_{-\beta}(\sigma)\nu(t)dt, \qquad (9)$$

где $\gamma_2 = 2\Gamma(2-2\alpha)/\Gamma^2(3/2-\alpha)$, $\sigma = \lambda\sqrt{(\eta-t)(t-\xi)}$,

$$A_\alpha^-(\tau,\lambda) = \gamma_1(\eta-\xi)^{-1-2\beta}\int_\xi^\eta (\eta-t)^\beta(t-\xi)^\beta \overline{J}_\beta(\sigma)\tau(t)dt - $$

$$-\frac{\gamma_1(\eta-\xi)^{-1-2\beta}}{2(1+2\beta)(\beta+1)}\int_\xi^\eta(\eta-t)^{\beta+1}(t-\xi)^{\beta+1}\overline{J}_{\beta+1}(\sigma)\left[\lambda^2\tau(t) - \tau''(t)\right]dt.$$

**Определение 2.** *Функция $u(x,y)$, определяемая в области $\Omega_2$ формулой (9), называется решением уравнения $L_{\alpha,\lambda}(u)=0$ из класса $R_{0p}^\lambda$ при $-1/2 < \alpha < 0$, если функция $\tau(x)$ представима в виде*

$$\tau(x) = sign(x-p)\int_p^x |x-t|^{-2\beta}\overline{I}_{-\beta}\left[\lambda(x-t)\right]T(t)dt,$$

*где $\nu(x), T(x) \in C[0,1]\cap C^1(0,1)$ и $\nu'(x), T'(x) \in L(0,1)$.*

Согласно определению 2, функция $u(x,y)$, определенная в области $\Omega_2$ в виде (9), называется решением уравнения $L_{\alpha,\lambda}(0)=0$ из класса $R_{00}^\lambda$, если функция $\tau(x)$ представима в виде

$$\tau(x) = \int_0^x (x-s)^{-2\beta}\overline{I}_{-\beta}\left[\lambda(x-s)\right]T(s)ds, \qquad (10)$$

а $\nu(x), T(x) \in C[0,1]\cap C^1(0,1)$ и $\nu'(x), T'(x) \in L(0,1)$.

Из (10) вытекает, $\tau(x) \in C^3[0,1]$ и

$$\tau(0) = 0,\ \tau'(0) = 0. \qquad (11)$$

Решения задачи {(1), (7), (8)} из класса $R_{00}^\lambda$ в области $\Omega_2$ имеет вид

$$u(x,y) = \int_0^\xi (\eta-s)^{-\beta}(\xi-s)^{-\beta}\overline{I}_{-\beta}\left[\lambda\sqrt{(\eta-s)(\xi-s)}\right]T(s)ds +$$

$$+ \int_\xi^\eta (\eta-s)^{-\beta}(s-\xi)^{-\beta}\overline{J}_{-\beta}\left[\lambda\sqrt{(\eta-s)(s-\xi)}\right]N(s)ds, \qquad (12)$$

где $N(s) = (2cos\pi\beta)^{-1}T(s) - 4^{2\beta-1}\gamma_2 \nu(s)$.

Подчиняя решение (12) краевому условию (3), получаем уравнение относительно $N(\eta)$:

$$\int_0^\eta (\eta-s)^{-\beta} s^{-\beta} \overline{J}_{-\beta}\left[\lambda\sqrt{(\eta-s)s}\right]N(s)ds = \psi\left(\frac{\eta}{2}\right), \quad 0 \le \eta \le 1.$$

Последнее уравнение, в результате применения равенства (6), можно привести к виду, удобному для дальнейшего исследования:

$$D_{0\eta}^{\beta-1}\left\{B_{0\eta}^{1,\lambda i}\left[\eta^{-\beta}N(\eta)\right]\right\} = \frac{1}{\Gamma(-\beta+1)}\psi\left(\frac{\eta}{2}\right), \quad 0 \le \eta \le 1. \qquad (13)$$

Применяя к обеим частям уравнения (13) последовательно операторы $D_{0\eta}^{1-\beta}$, $A_{0\eta}^{1,\lambda i}$ и учитывая равенства $D_{0\eta}^{1-\beta}D_{0\eta}^{\beta-1}f(\eta) = f(\eta)$, $A_{0\eta}^{1,\lambda i}B_{0\eta}^{1,\lambda i}g(\eta) = g(\eta)$, а также структуры функции $N(\eta)$, получаем

$$T(x) = \gamma_3 \nu(x) + \frac{2cos\pi\beta x^\beta}{\Gamma(1-\beta)} A_{0x}^{1,i\lambda}\left[D_{0x}^{1-\beta}\psi\left(\frac{x}{2}\right)\right].$$

Подставляя это значение $T(x)$ в (10), находим соотношение между $\tau(x)$ и $\nu(x)$ на отрезке $[0,1]$, получаемое из области $\Omega_2$:

$$\tau(x) = \gamma_3 \int_0^x (x-s)^{-2\beta}\overline{I}_{-\beta}\left[\lambda(x-s)\right]\nu(s)ds +$$

$$+ \frac{2\cos\pi\beta}{\Gamma(1-\beta)}\int_0^x (x-s)^{-2\beta}\overline{I}_{-\beta}\left[\lambda(x-s)\right]s^\beta A_{0s}^{1,\lambda}\left[D_{0\zeta}^{1-\beta}\psi\left(\frac{s}{2}\right)\right]ds, \quad 0 \leqslant x \leqslant 1, \qquad (14)$$

где $\gamma_3 = 2\cdot 4^{2\beta-1}\gamma_2 cos\pi\beta$.

Согласно условиям задачи $T_0$, в уравнении (1) и в условиях (2) можно перейти к пределу при $y \to +0$ (например, см. [6], [7]). В результате получим следующие соотношения:

$$\tau''(x) - \Gamma(1+\delta)\nu(x) - \lambda^2 \tau(x) = 0, \quad 0 \leqslant x \leqslant 1, \qquad (15)$$

$$\tau(0) = \lim_{y \to +0} y^{1-\delta} \varphi_1(y) = a, \ \tau(1) = \lim_{y \to +0} y^{1-\delta} \varphi_2(y) = b. \qquad (16)$$

Если считать, что $\nu(x)$ - известная функция, то при $\lambda \in R$ или $\lambda i \in R$, $\lambda i \neq \pi m, m \in Z$ задача $\{(15), (16)\}$ имеет единственное решение [8]

$$\tau(x) = a + x(b-a) +$$

$$+ \lambda^2 \int_0^1 G(x,t;\lambda)\bigl[a + t(b-a)\bigr]dt + \Gamma(1+\delta)\int_0^1 G(x,t;\lambda)\nu(t)dt, \qquad (17)$$

где $G(x,t;\lambda)$ - функция Грина задачи $\{(15), (16)\}$

$$G(x,t;\lambda) = \begin{cases} \dfrac{sh\lambda(x-t) sh\lambda t}{\lambda sh\lambda}, 0 \leq x \leq t, \\ \dfrac{sh\lambda x\, sh\lambda(t-1)}{\lambda sh\lambda}, t \leq x \leq 1. \end{cases}$$

Из формулы (17), в силу равенства (11), вытекают следующие равенства $a = 0, b = 0$ и

$$\tau(x) = \Gamma(1+\delta)\int_0^1 G(x,t;\lambda)\nu(t)dt. \qquad (18)$$

(18) является основным соотношением между $\tau(x)$ и $\nu(x)$ на отрезке $[0,1]$, получаемое из области $\Omega_1$.

Теперь из соотношений (18) и (14) найдем неизвестные функции $\tau(x)$ и $\nu(x)$. С этой целью, из (18) и (14) исключим функцию $\tau(x)$:

$$\Gamma(1+\delta)\int_0^1 G(x,t;\lambda)\nu(t)dt = \gamma_3 \int_0^x (x-\zeta)^{-2\beta} \overline{I}_{-\beta}\bigl[\lambda(x-\zeta)\bigr]\nu(\zeta)d\zeta +$$

$$+ \frac{2\cos\pi\beta}{\Gamma(1-\beta)}\int_0^x (x-s)^{-2\beta} \overline{I}_{-\beta}\bigl[\lambda(x-s)\bigr]s^\beta A_{0s}^{1,\lambda}\left[D_{0s}^{1-\beta}\psi\left(\frac{s}{2}\right)\right]ds, \ 0 \leqslant x \leqslant 1, \quad (19)$$

Продифференцируем это равенство дважды по $x$. Затем, от полученного равенства почленно вычтем равенство (19). В результате, имеем интегральное уравнение относительно $\nu(x)$:

$$\nu(x) - \frac{2\beta(2\beta+1)\gamma_3}{\Gamma(1+\delta)}\int_0^x (x-s)^{-2\beta-2} \overline{I}_{-\beta-1}\bigl[\lambda(x-s)\bigr]\nu(s)ds = Q(x), 0 \leq x \leq 1, \quad (20)$$

где

$$Q(x) = \frac{4\beta(2\beta+1)\cos\pi\beta}{\Gamma(1+\delta)\Gamma(1-\beta)} \int_0^x (x-t)^{-2\beta-2} \overline{I}_{-\beta-1}\left[\lambda(x-t)\right] \Phi(t) dt,$$

$$\Phi(t) = t^\beta A_{0t}^{1,\lambda}\left[D_{0t}^{1-\beta}\psi\left(\frac{t}{2}\right)\right].$$

Так как $\alpha \in (-1/2, 0)$, то $\beta = \alpha - 1/2 \in (-1, -1/2)$ и $-2\beta - 2 \in (-1, 0)$. Поэтому ядро интегрального уравнения (20) имеет слабую особенность.

Пусть выполнены следующие условия:

$$\psi^{(m)}(0) = 0,\, m = 0, 1, 2,\, \psi'''(s/2) = s^p \psi_0(s),\, p > -2 - 2\beta,\, \psi_0(s) \in C[0,1]. \quad (21)$$

Докажем, что $Q(x) \in C[0,1] \cap C^1(0,1)$ и $Q'(x) \in L(0,1)$. С учетом (21) и (4), имеем

$$\Phi(t) = \frac{t^{2+2\beta}}{8\Gamma(2+\beta)} \int_0^1 (1-z)^{1+\beta} \psi'''\left(\frac{zt}{2}\right) dz -$$

$$- \frac{\lambda^2 t^{4+2\beta}}{32\Gamma(2+\beta)} \int_0^1 s^{3+\beta} \overline{J}_1\left[\lambda t\sqrt{s(1-s)}\right] ds \int_0^1 (1-z)^{1+\beta} \psi'''\left(\frac{zts}{2}\right) dz,$$

$$\Phi'(t) = \frac{\beta t^{1+2\beta}}{8\Gamma(2+\beta)} \int_0^1 (1-z)^{1+\beta} \psi'''\left(\frac{tz}{2}\right) dz + \frac{t^{1+2\beta}}{8\Gamma(1+\beta)} \int_0^1 (1-z)^\beta \psi'''\left(\frac{zt}{2}\right) dz -$$

$$- \frac{\beta\lambda^2 t^{3+2\beta}}{32\Gamma(2+\beta)} \int_0^1 s^{3+\beta} \overline{J}_1\left[\lambda t\sqrt{s(1-s)}\right] ds \int_0^1 (1-z)^{1+\beta} \psi'''\left(\frac{zst}{2}\right) dz +$$

$$+ \frac{\lambda^4 t^{5+2\beta}}{128\Gamma(2+\beta)} \int_0^1 s^{3+\beta}(2-s) \overline{J}_2\left[\lambda t\sqrt{s(1-s)}\right] ds \int_0^1 (1-z)^{1+\beta} \psi'''\left(\frac{zts}{2}\right) dz. \quad (22)$$

Отсюда, в силу $\psi'''(s/2) = s^p \psi_0(s), p > -2 - 2\beta, \psi_0(s) \in C[0,1]$, следует, что $\Phi(\zeta) \in C[0,1]$, поэтому $Q(x) \in C[0,1]$. Теперь вычисляем $Q'(x)$:

$$Q'(x) = -\frac{4\beta \cos\pi\beta \lambda^2}{\Gamma(1-\beta)} \int_0^x (x-t)^{-2\beta-1} \overline{I}_{-\beta-1}\left[\lambda(x-t)\right] \Phi(t) dt +$$

$$+ \frac{4\beta(2\beta+1)\cos\pi\beta}{\Gamma(1-\beta)} \int_0^x (x-t)^{-2\beta-2} \overline{I}_{-\beta-1}\left[\lambda(x-t)\right] \Phi'(t) dt -$$

$$- \frac{4\beta \cos\pi\beta \lambda}{\Gamma(1-\beta)} \int_0^x (x-t)^{-2\beta-1} \overline{I}'_{-\beta-1}\left[\lambda(x-t)\right] \Phi'(t) dt.$$

Отсюда, согласно (21) и (22) следует, что $Q'(x) \in C(0,1) \cap L(0,1)$. Следовательно, $Q(x) \in C[0,1] \cap C^1(0,1)$ и $Q'(x) \in L(0,1)$.

В силу свойств ядра и правой части интегрального уравнения (20), согласно теории интегральных уравнений Вольтерра второго рода [9], оно имеет единственное решение.

После нахождения функции $\nu(x)$ из (20), функция $\tau(x)$ находится по формуле (18). После этого решение задачи $T_0$ в области $\Omega_2$ определяется по формуле (12), а в области $\Omega_1$ - как решение первой краевой задачи для уравнения $L_1(u) = 0$ с условиями (2) и $\lim\limits_{y \to +0} y^{1-\delta} u(x,y) = \tau(x)$, $0 \leqslant x \leqslant 1$, определяется по формуле [10]

$$u(x,y) = \int_0^1 \tau(t) G(x,y;t,0) dt + \int_0^y \varphi_1(s) G_t(x,y;0,s) ds - \int_0^y \varphi_2(s) G_t(x,y;1,s) ds,$$

где

$$G(x,y;t,s) = \sum_{m=-\infty}^{+\infty} \left[ \Gamma(x-t+2m, y-s) - \Gamma(x+t+2m, y-s) \right],$$

$$\Gamma(x,y) = \frac{1}{2y} \int_{|x|}^\infty e_{1,\delta/2}^{1,0}\left( -\frac{\xi}{y^{\delta/2}} \right) J_0\left( \lambda \sqrt{\xi^2 - x^2} \right) d\xi, \quad e_{\alpha,\beta}^{\mu,\delta}(z) = \sum_{k=0}^\infty \frac{z^k}{\Gamma(\mu + \alpha k)\Gamma(\delta - \beta k)}, \alpha > \beta.$$

Таким образом, доказана следующая основная

**Теорема 2.** *Если $\lambda$ – действительное число или чисто мнимое число, отличное от $i\pi n$, $n \in Z$, а заданные функции удовлетворяют условиям (21) и $y^{1-\delta} \varphi_1(y), y^{1-\delta} \varphi_2(y) \in C[0,1]$, $\lim\limits_{y \to +0} y^{1-\delta} \varphi_1(y) = 0$, $\lim\limits_{y \to +0} y^{1-\delta} \varphi_2(y) = 0$, то задача $T_0$ имеет единственное решение.*